\documentclass[11pt]{article}

\usepackage{amsfonts}
\usepackage{graphics}
\usepackage{indentfirst}
\usepackage{color}
\usepackage{cite}
\usepackage{float}
\usepackage{caption}
\usepackage{latexsym}
\usepackage[paper=a4paper, left=1.6cm, right=1.6cm, top=1.8cm, bottom=1.6cm, headheight=5.5pt, footskip=0.8cm, footnotesep=0.8cm, centering, includefoot]{geometry}
\usepackage{amsmath}
\allowdisplaybreaks
\usepackage{amssymb}
\usepackage[colorlinks, linkcolor=red]{hyperref}
\hypersetup{urlcolor=red, citecolor=red}
\usepackage[dvips]{epsfig}
\usepackage{amscd}
\usepackage{amsthm}
\usepackage{mathrsfs}
\usepackage{verbatim}

\newtheorem{theorem}{Theorem}[section]
\newtheorem{remark}{Remark}[section]
\newtheorem{definition}{Definition}[section]
\newtheorem{lemma}{Lemma}[section]
\newtheorem{corollary}{Corollary}[section]
\newtheorem{proposition}{Proposition}[section]

\newcommand{\cH}{\mathscr H}

\title{Scaling-invariant smallness and global strong solutions
	to the three-dimensional isentropic compressible Navier--Stokes equations with a non-vacuum far-field density
	\thanks{This research was partially supported by National Natural Science Foundation of China (No. 12371227) and Fundamental Research Funds for the Central Universities (No. SWU--KU24001).}
}

\author{Lin Xu,\ Xin Zhong{\thanks{Corresponding author. E-mail addresses: mathxu@email.swu.edu.cn (L. Xu), xzhong1014@amss.ac.cn (X. Zhong).}}
	\date{}\\
	\footnotesize School of Mathematics and Statistics, Southwest University, Chongqing 400715, P. R. China}
\begin{document}
	\maketitle
	
	\begin{abstract}
		We consider the Cauchy problem for the three-dimensional
		isentropic compressible Navier--Stokes equations with a
		positive far-field density $\rho_\infty$.
		We identify the following scaling-invariant quantity
		associated with the initial data:
	\begin{align*}
		\rho_\infty^{\frac{3-\gamma}{2}}
		\left(
		\|\nabla u_0\|_{L^2}^2
		+P(\rho_\infty)
		\Big(
		\frac12\int_{\mathbb{R}^3}\rho_0|u_0|^2\,dx
		+\int_{\mathbb{R}^3}G(\rho_0)\,dx
		\Big)
		\right),
	\end{align*}
		where $G$ denotes the relative potential energy density
		with respect to $\rho_\infty$.
		Provided that the scaling-invariant quantity is sufficiently small,
		we establish the global existence and uniqueness of strong solutions, allowing the initial density to contain interior vacuum and exhibit large oscillations.
		The smallness criterion depends explicitly on the far-field density through the prefactor $\rho_\infty^{(3-\gamma)/2}$ and the background pressure $P(\rho_\infty)$.
		Our result provides a non-vacuum counterpart to Wen's
		scaling-invariant global existence result for the
		far-field vacuum case
		(\textit{Adv. Math.} \textbf{482} (2025), Paper No. 110628). 
	\end{abstract}



\setcounter{tocdepth}{2}
\tableofcontents

\section{Introduction and main result}

We consider the Cauchy problem for the three-dimensional
isentropic compressible Navier--Stokes equations in
$\mathbb R^3$:
\begin{align}\label{1}
	\begin{cases}
		\rho_t+\operatorname{div}(\rho u)=0,\\
		(\rho u)_t+\operatorname{div}(\rho u\otimes u)
		-\mu\Delta u-(\mu+\lambda)\nabla\operatorname{div}u
		+\nabla P(\rho)=0,
	\end{cases}
\end{align}
where $\rho=\rho(x,t)\geq0$ and
$u=u(x,t)\in\mathbb R^3$ denote the fluid density and
velocity, respectively.
The pressure is given by the polytropic law
\begin{align*}
	P(\rho)=A\rho^\gamma,\qquad A>0,\quad\gamma>1,
\end{align*}
while the constant viscosity coefficients $\mu$ and
$\lambda$ satisfy the physical restrictions
\begin{align*}
	\mu>0,\qquad 2\mu+3\lambda\geq0.
\end{align*}
We prescribe the initial data
\begin{align}\label{2}
	(\rho,~ \rho u)|_{t=0}=(\rho_0, ~\rho_0u_0),
\end{align}
and impose the non-vacuum far-field condition
\begin{align}\label{3}
	\lim_{|x|\to\infty}(\rho(x,t),~u(x,t))
	=(\rho_\infty,~0), 
\end{align}
where $\rho_\infty>0$ is a given constant representing the far-field density.

Developing a global theory of smooth solutions to the
multidimensional compressible Navier--Stokes equations
for general large initial data remains an
open problem in the mathematical theory of viscous
compressible fluids.
One classical approach, initiated by Matsumura and
Nishida \cite{MN1980}, is based on small perturbations
of a   positive constant equilibrium state.
Another important direction concerns global solutions
whose initial density may contain vacuum and exhibit large
oscillations.
In this setting, Huang, Li, and Xin \cite{HLX2012}
established the global existence and uniqueness of
classical solutions under a suitable smallness condition
on the initial energy, allowing either vacuum or non-vacuum state at infinity.
Further developments concerning global solutions with
vacuum or certain classes of large initial data can be
found in \cite{LiXin2019,HHPZ2024,HL2022}.
These results demonstrated the
importance of suitable smallness assumptions or
additional structural properties in establishing
global well-posedness.

A different approach to the global well-posedness problem
is based on critical function spaces.
Danchin \cite{Danchin2000} established global
well-posedness in critical Besov spaces for small
perturbations of a positive constant state.
This theory was subsequently extended to critical
$L^p$-type Besov spaces by Charve and Danchin \cite{CD2010}
and Chen, Miao, and Zhang \cite{CMZ2010},  they obtained
global solutions with highly oscillating initial
velocities.
Haspot further enlarged the admissible critical spaces in
\cite{H2011a} and established global strong solutions
for a   class  of large initial data under
additional regularity assumptions in \cite{H2011b}.
For further developments in Besov spaces, see \cite{FZZ2018,HHW2019,ZLZ2020,WWZ2014}.
These results demonstrate the effectiveness of
critical-space methods in treating highly oscillatory
velocity or large initial data with additional assumptions.
In the critical-space framework, smallness is typically imposed on perturbations around positive constant states.
This  framework is  different from
energy-based approaches that allow  density vacuum
and large  oscillations.

\textit{A natural question is whether one can establish
	global strong solutions under an energy-based
	smallness condition that is invariant under the
	intrinsic scaling of system \eqref{1}.}
For the polytropic pressure law, such a scaling
was identified by Lei and Xin \cite{LX2019}.
More precisely, the system \eqref{1} is invariant under
the transformation (details see  \cite[Appendix]{XuZhong2026}):
\begin{align}\label{scaling}
	\rho^\ell(x,t)
	\triangleq\ell^{\frac{2}{\gamma+1}}
	\rho\left(
	\ell x,
	\ell^{\frac{2\gamma}{\gamma+1}}t
	\right),\quad
	u^\ell(x,t)
	\triangleq\ell^{\frac{\gamma-1}{\gamma+1}}
	u\left(
	\ell x,
	\ell^{\frac{2\gamma}{\gamma+1}}t
	\right),
	\qquad \ell>0.
\end{align}
An essential feature of this transformation is
that the pressure gradient scales at the same
rate as the inertial and viscous terms.
This differs from the usual parabolic scaling
\begin{align*}
	\rho^k(x,t)
	\triangleq\rho(k x,k^2t),\quad
	u^k(x,t)
	\triangleq k u(k x,k^2t),\quad k>0,
\end{align*}
which is commonly used to identify critical
regularity in the critical framework but
does not preserve the polytropic pressure law $P=A\rho^{\gamma}$ scaling invariant.
Moreover, the initial energy is not invariant under the scaling \eqref{scaling} and its smallness alone does not yield a scaling-invariant quantity. This suggests combining the energy with suitable norms of the initial data to form a scaling-invariant quantity. For the heat-conducting compressible Navier--Stokes
system, Li \cite{Li2020} established
global strong solutions under a smallness condition invariant
under the natural scaling of that system.
For system \eqref{1} with far-field vacuum, Wen
\cite{Wen2025} established global strong solutions under
an explicit smallness condition invariant under
\eqref{scaling}.
In our recent work \cite{XuZhong2026}, we
investigated the role of pressure in determining
the critical scaling structure and developed
a corresponding global low-regularity theory.
For further related results, see \cite{XXZ2026,XZ20261,XZ202611}.
These developments demonstrate the effectiveness of scaling-invariant smallness conditions in establishing global strong solutions and complement the small-energy theory developed in \cite{HLX2012}.
 
\subsection{Motivation}
The preceding developments lead us to investigate
how a scaling-invariant smallness condition can be
formulated for system \eqref{1}
with a nonzero far-field density. It is worth emphasizing that the non-vacuum far-field problem exhibits essential differences from the far-field vacuum case. Under the intrinsic scaling \eqref{scaling}, the positive far-field density $\rho_\infty$ changes accordingly, while the estimates for the pressure perturbation and the relative potential energy must account for the contribution of the nonzero background pressure. These features require a separate analysis of the non-vacuum setting, since the arguments developed for the far-field vacuum case cannot be directly applied.

\begin{itemize}
	\item Indeed, when $\rho_\infty=0$, the far-field vacuum state is preserved under the scaling \eqref{scaling}. In contrast, when $\rho_\infty>0$, the far-field density transforms according to
	\begin{align*}
	\rho_\infty^\ell
	=\ell^{\frac{2}{\gamma+1}}\rho_\infty.
\end{align*}
	Consequently, the prescribed far-field density is not invariant under the scaling if $\rho_\infty$ is kept fixed. Nevertheless, since $\rho_\infty^\ell>0$ for every $\ell>0$, the scaling maps a Cauchy problem with a positive far-field density to another problem in the same class, with the background density transformed accordingly.
	In particular, setting
	$\bar\rho\triangleq\|\rho_0\|_{L^\infty}$, we observe that
	\begin{align}\label{ss}
		R^\ell
		\triangleq
		\frac{\bar\rho^\ell}{\rho_\infty^\ell}=
		\frac{
			\ell^{\frac{2}{\gamma+1}}\bar\rho
		}{
			\ell^{\frac{2}{\gamma+1}}\rho_\infty
		}=
		\frac{\bar\rho}{\rho_\infty}
		\triangleq R
	\end{align}		
	remains unchanged under the transformation \eqref{scaling}.
	This observation suggests that $\rho_\infty$ can serve as a natural
	normalization parameter in the construction
	of scaling-invariant initial quantities.
	
	\item  A positive far-field density also gives rise to
	a nonzero background pressure $P(\rho_\infty)$.
	Unlike the far-field vacuum case \cite{Wen2025,XuZhong2026,XZ20261}, the energy
	framework must be formulated in terms of
	perturbations relative to the positive background
	state, involving
	$P(\rho)-P(\rho_\infty)$ and the associated
	potential energy.
\end{itemize}	
These considerations motivate us to construct a
 scaling-invariant initial
quantity that combines the velocity,
the pressure perturbation, and the relative energy.
Under a suitable smallness assumption on
this quantity, we establish the global
existence and uniqueness of strong solutions
to the Cauchy problem \eqref{1}--\eqref{3}.

\subsection{Main result}
Before stating our main result, we introduce some notations used throughout
the paper.
The notation
$a\triangleq b$ means that $a$ is defined to be $b$.
For two nonnegative quantities $a$ and $b$, we write
$a\sim b$ if
$
cb\le a\le Cb$
for some positive constants $c$ and $C$.
Unless otherwise specified, $c$ and $C$ denote generic
positive constants that may vary from line to line.
Their dependence on the parameters will be specified
when relevant.  For an integer $k\geq0$ and
$1\leq q\leq\infty$, we define
\begin{gather*}
	D^{k,q}(\mathbb{R}^{3}) = \{ u \in L^1_{loc}(\mathbb{R}^{3}) \mid \|\nabla^k u\|_{L^q} < \infty \}, \ \|u\|_{D^{k,q}} = \|\nabla^k u\|_{L^{q}(\mathbb{R}^{3})},\ \|  u \|_{L^{q}}=\|u \|_{L^{q}(\mathbb{R}^{3})}, \\
	W^{k,q}(\mathbb{R}^{3}) = L^q(\mathbb{R}^{3}) \cap D^{k,q}(\mathbb{R}^{3}), \ D^k(\mathbb{R}^{3}) = D^{k,2}(\mathbb{R}^{3}), \ H^k(\mathbb{R}^{3}) = W^{k,2}(\mathbb{R}^{3}).
\end{gather*} We also use the shorthand
notation
\begin{align*}
	\int \cdot \, d x \triangleq \int_{\mathbb{R}^{3}} \cdot \, d x,\quad  \dot{u}\triangleq u_t+u\cdot\nabla u,\quad 	L^p\triangleq L^p(\mathbb{R}^{3}),~~ W^{k,p}\triangleq W^{k,p}(\mathbb{R}^{3}),~~ H^k\triangleq W^{k,2}(\mathbb{R}^{3}).
\end{align*}

We set
\begin{align*}
P_\infty
\triangleq P(\rho_\infty)
=A\rho_\infty^\gamma,
\qquad
p=p(t)\triangleq P(\rho(t))-P_\infty,
\qquad
\nu\triangleq2\mu+\lambda>0.
\end{align*}
The effective viscous flux and vorticity are defined by
\begin{align}\label{eqf}
	F\triangleq\nu\operatorname{div}u-p,
	\qquad
	\omega\triangleq\operatorname{curl}u
	=\nabla\times u.
\end{align}
We define the relative energy functional $E(t)$ and the
higher-order quantity $D(t)$ by
\begin{align*}
	E(t)\triangleq\frac12\int \rho|u|^2(t)\, dx+\int G(\rho(t))\, dx,\quad
	D(t)\triangleq\|\nabla u(t)\|_{L^2}^2+\|p(t)\|_{L^2}^2,
\end{align*}
where $G$ denotes the relative potential-energy density,
which is given by
\begin{align*}
	G(\rho)\triangleq\rho\int_{\rho_\infty}^{\rho}\frac{P(s)-P_\infty}{s^2}\, ds
	=\frac{A}{\gamma-1}\left(\rho^\gamma-\gamma\rho_\infty^{\gamma-1}\rho+(\gamma-1)\rho_\infty^\gamma\right).
\end{align*}

The strong solutions considered in this paper are defined
as follows.
	\begin{definition}\label{def}
		Let $T>0$ and $q\in(3,6)$. Assume that the initial data
		$(\rho_0\geq 0,u_0)$ satisfy
		\begin{align}\label{sou}
			\left( \rho_{0}-\rho_{\infty}\right) \in   H^{1}\cap W^{1,q},\qquad u_{0}\in D^{1},\qquad \rho_{0}u_{0}\in L^{2}.
		\end{align}
		A pair $(\rho,u)$ is called a strong solution to system  \eqref{1}  in $\mathbb{R}^{3}\times[0,T]$ if it satisfies the following regularity properties:
		\begin{align}\label{defs}
			\begin{cases}
				(\rho-\rho_{\infty})\in C\big([0,T];L^{2}\big)
				\cap L^{\infty}\big(0,T;H^{1}\cap W^{1,q}\big),
				\quad
				\rho_{t}\in L^{\infty}\big(0,T;L^{2}\big), \\
				\rho u\in C\big([0,T];L^{2}\big),\quad
				u\in L^{\infty}\big(0,T;D^{1}\big)
				\cap L^{2}\big(0,T;D^{2}\big),
				\quad
				\sqrt{\rho} u_{t}\in L^{2}\big(0,T;L^{2}\big), \\
				\sqrt{t} u\in L^{\infty}\big(0,T;D^{2}\big)
				\cap L^{2}\big(0,T;D^{2,q}\big),
				\quad
				\sqrt{t} u_{t}\in L^{2}\big(0,T;D^{1}\big).
			\end{cases}
		\end{align}
		Moreover, $(\rho,u)$ satisfies \eqref{1}  a.e. in $\mathbb{R}^{3}\times(0,T]$. In particular, the strong solution  $(\rho, u )$  of \eqref{1}  is called the global strong solution, if the strong solution satisfies \eqref{defs} for any  $T>0$, and satisfies \eqref{1} a.e. in  $\mathbb{R}^{3} \times(0, T)$.
	\end{definition}

	Our main result is the following.
	
	\begin{theorem}\label{th1}
		For $R= \bar\rho/\rho_\infty$, $\bar\rho=\|\rho_0\|_{L^\infty}$. Assume that the initial data $(\rho_0,u_0)$ satisfy
		\eqref{sou}. Then there exists a positive constant
		\begin{align*}
		\varepsilon_0
		=\varepsilon_0(R,\gamma,A,\mu,\lambda)
		\end{align*}
		such that, if
		\begin{align}\label{sca}
			\mathcal S_{0}
			\triangleq
			\rho_\infty^{\frac{3-\gamma}{2}}
			\bigl[\|\nabla u_{0}\|_{L^2}^2+\|P(\rho_0)-P_\infty\|_{L^2}^2+P_\infty E(0)\bigr]
			\leq\varepsilon_0,
		\end{align}
		the Cauchy problem \eqref{1}--\eqref{3}
		admits a unique global strong solution $(\rho,u)$.
	\end{theorem}

	An important observation is that the initial pressure
	perturbation $\|P(\rho_0)-P_\infty\|_{L^2}^2$ appearing in $\mathcal S_{0}$ can be controlled by
	the relative potential energy, provided that the
	initial density is bounded from above.
	This leads to the following equivalent formulation
	of the smallness condition.
	
	\begin{corollary} 
		There exists a constant $C_{\gamma,R}>0$,
		depending only on $\gamma$ and $R$, such that
		\begin{align*}
			\widehat{\mathcal S_{0}}
			\leq \mathcal S_{0}
			\leq(1+C_{\gamma,R})\widehat{\mathcal S_{0}},
		\end{align*}
		where
		\begin{align*}
		\widehat{\mathcal S_{0}}
			\triangleq
			\rho_\infty^{\frac{3-\gamma}{2}}
			\left(
			\|\nabla u_0\|_{L^2}^2
			+P_\infty E(0)
			\right).
		\end{align*}
		Then, under the assumptions of
		Theorem \ref{th1}, the conclusion of the theorem
		remains valid if the condition \eqref{sca} is replaced by
		\begin{align}\label{sca2}
			\widehat{\mathcal S_{0}}
			\leq\frac{\varepsilon_0}{1+C_{\gamma,R}}.
		\end{align}
	\end{corollary}
	
	\begin{proof}
		Set $s_{0}\triangleq \rho_0/\rho_\infty$. Since
		$0\leq\rho_0\leq\bar\rho$, we have $0\leq s_{0}\leq R$.
		Recalling the definitions of the relative potential
		energy and the pressure perturbation, we obtain
		\begin{align*}
		G(\rho_0)
		=\frac{P_\infty}{\gamma-1}
		\left(s_{0}^\gamma-\gamma s_{0}+\gamma-1\right)\quad \text{and}\quad
		p(0)=P(\rho_0)-P_\infty
		=P_\infty(s_{0}^\gamma-1).
		\end{align*}
		Since $\gamma>1$, the strict convexity of
		$s_{0}\mapsto s_{0}^\gamma$ implies that
		\begin{align*}
		s_{0}^\gamma-\gamma s_{0}+\gamma-1>0
		\qquad\text{for }s_{0}\geq0,\ s_{0}\neq1.
		\end{align*}
		Moreover, a Taylor expansion around $s_{0}=1$ yields
		\begin{align*}
		\lim_{s_{0}\to1}
		\frac{(\gamma-1)(s_{0}^\gamma-1)^2}
		{s_{0}^\gamma-\gamma s_{0}+\gamma-1}
		=2\gamma.
		\end{align*}
		Consequently, by continuity on the interval
		$[0,R]$, there exists a constant
		$C_{\gamma,R}>0$ such that
		\begin{align*}
		(s_{0}^\gamma-1)^2
		\leq
		\frac{C_{\gamma,R}}{\gamma-1}
		\left(s_{0}^\gamma-\gamma s_{0}+\gamma-1\right),\quad \text{for all}\quad s_{0}\in[0,R].
		\end{align*}
		
		It follows that
		\begin{align*}
			\|p(0)\|_{L^2}^2
			\leq C_{\gamma,R}P_\infty
			\int G(\rho_0)\,dx\leq C_{\gamma,R}P_\infty E(0).
		\end{align*}
		Recalling that
		$
		\|p(0)\|_{L^2}^2
		=\|P(\rho_0)-P_\infty\|_{L^2}^2$,
		we immediately obtain
		\begin{align*}
		\widehat{\mathcal S_{0}}
		\leq \mathcal S_{0}
		\leq(1+C_{\gamma,R})\widehat{\mathcal  S_{0}}.
		\end{align*}
		In particular, \eqref{sca2} implies
		$\mathcal S_{0}\leq\varepsilon_0$.
		The desired conclusion therefore follows
		from Theorem \ref{th1}.
	\end{proof}
	\begin{remark}[Scaling invariance]
 It is worth mentioning that the quantity appeared in   \eqref{sca} is invariant under the scaling \eqref{scaling}. In particular, the  threshold $\varepsilon_0$ is permitted to depend on $R$, which is itself invariant under the scaling. Under \eqref{scaling}, a direct calculation gives
		\begin{align*}
			E^\ell(t)=\ell^{-\frac{\gamma+3}{\gamma+1}}
			E\left(\ell^{\frac{2\gamma}{\gamma+1}}t\right),\quad
			D^\ell(t)=\ell^{\frac{\gamma-3}{\gamma+1}}
			D\left(\ell^{\frac{2\gamma}{\gamma+1}}t\right),\quad
			P_\infty^\ell=\ell^{\frac{2\gamma}{\gamma+1}}P_\infty,
		\end{align*} 
		we have $\mathcal{S}_0^\ell=\mathcal{S}_0$ \textnormal{(details see Lemma \ref{check})}.
	\end{remark}
	
	\begin{remark}
	 Since 	$
			\rho_{\infty}u_{0}=(\rho_{\infty}-\rho_{0})u_{0}+\rho_{0}u_{0}$,
		we have
		\begin{align*}
			\rho_{\infty}^{2}\|u_{0}\|_{L^{2}}^{2}
			\leq
			C\|(\rho_{\infty}-\rho_{0})u_{0}\|_{L^{2}}^{2}
			+C\|\rho_{0}u_{0}\|_{L^{2}}^{2}
			C\|\rho_{\infty}-\rho_{0}\|_{L^{3}}^{2}
			\|\nabla u_{0}\|_{L^{2}}^{2}
			+C\|\rho_{0}u_{0}\|_{L^{2}}^{2}.
		\end{align*}
		Moreover, since
		$\rho_{0}-\rho_{\infty}\in H^{1}\cap W^{1,q}$, $u_{0}\in D^{1}$, we have
		$\rho_{0}\in L^{\infty}$ and $u_{0}\in L^{2}$. Hence
		\begin{align*}
			\|\sqrt{\rho_{0}}u_{0}\|_{L^{2}}^{2}
			\leq
			\|\rho_{0}\|_{L^{\infty}}\|u_{0}\|_{L^{2}}^{2}
			\leq C,
		\end{align*}
		the initial energy $E_{0}$ is well-defined.
	\end{remark}
	
	\begin{remark}[Large initial energy along scaling orbits]
		An important feature of the scaling-invariant smallness
		condition is that it allows large initial
		energy along scaling-related families of Cauchy problems.
		
		Indeed, under the scaling transformation \eqref{scaling},
		the initial energy satisfies \textnormal{(see \eqref{eqE})}
		\begin{align*}
			E^\ell(0)
			=
			\ell^{-\frac{\gamma+3}{\gamma+1}}E(0).
		\end{align*}
		Consequently, for any initial data with $E(0)>0$,
		we have
		\begin{align*}
			E^\ell(0)\to\infty
			\qquad\text{as }\ell\to0^+,
		\end{align*}
		whereas the scaling-invariant quantity remains unchanged:
		\begin{align*}
			\mathcal S_0^\ell=\mathcal S_0.
		\end{align*}
		This illustrates the distinction between
		smallness of the initial energy and smallness of a
		quantity adapted to the intrinsic scaling of the system.
	\end{remark}
 
\begin{remark}[Comparison with the far-field vacuum case]
	Wen \cite{Wen2025} established the global
	well-posedness of system \eqref{1} with far-field vacuum
	under a suitable smallness condition on a
	scaling-invariant initial quantity
	\begin{align*}
		\bar{\rho}^{3}
		\left( 
		\frac{1}{2}\|\sqrt{\rho_{0}}u_{0}\|_{L^{2}}^{2}
		+\frac{\|P(\rho_{0})\|_{L^{1}}}{\gamma-1}
		\right) 
		\left(
		\|\nabla u_{0}\|_{L^{2}}^{2}
		+\|P(\rho_{0})\|_{L^{2}}^{2}
		\right)
		\left[ 
		1+\bar{\rho}^{3+\gamma}
		\Big(
		\frac{1}{2}\|\sqrt{\rho_{0}}u_{0}\|_{L^{2}}^{2}
		+\frac{1}{\gamma-1}\|P(\rho_{0})\|_{L^{1}}
		\Big)^{2}
		\right] .
	\end{align*}
	 Theorem \ref{th1} provides a non-vacuum
	counterpart to this result.
	A principal difference lies in the treatment
	of the positive equilibrium state $\rho_\infty$ and the pressure perturbation $p=P(\rho)-P_\infty$ enters
	the energy estimates. In particular, the pressure perturbation
	satisfies
	\begin{align*}
	p_t+u\cdot\nabla p
	+\gamma p\,\operatorname{div}u
	=
	-\gamma P_\infty\operatorname{div}u.
     \end{align*}
	The term on the right-hand side reflects
	the contribution of the nonzero background
	pressure.
	Our estimates account for this contribution
	through a modified energy functional
	involving $P_\infty E(t)$ \textnormal{(see \eqref{hh})}.
\end{remark}

\begin{remark}[Comparison with the classical small-energy theory]
	Huang, Li, and Xin \cite{HLX2012} established
	the global existence of classical solutions
	with small initial energy.
	Our result is formulated in a different
	smallness framework. Instead of imposing
	smallness solely on the relative energy,
	we  provide a scaling-invariant
	energy framework adapted to the nonzero
	far-field state.
\end{remark}

\begin{remark}[Dependence on the viscosity coefficients] 
The smallness threshold $\varepsilon_0$ in Theorem \ref{th1} depends on the viscosity coefficients $\mu$ and $\lambda$. It is therefore natural to ask whether one can establish estimates that are uniform with respect to these coefficients and use them to investigate the corresponding singular limits. We leave this question for future work.
\end{remark}

\begin{remark}[The two-dimensional case]
	It is worth pointing out that the two-dimensional and
	three-dimensional problems exhibit essential differences,
	particularly in their scaling properties and the
	associated a priori estimates. This raises a natural and
	important question of whether a scaling-invariant
	smallness framework can also be established for the
	two-dimensional Cauchy problem. In this direction, we
	have recently completed a study of the two-dimensional
	isentropic compressible Navier--Stokes equations in both
	the far-field vacuum and non-vacuum settings. The
	corresponding results will be presented in a forthcoming
	work.
\end{remark}

\subsection{Strategy of the proof}

We briefly outline the main difficulties and key ideas in the
proof of Theorem~\ref{th1}. The proof relies on a bootstrap
argument; see Proposition \ref{pro}.
Compared with the far-field vacuum case, the presence
of a positive constant state $\rho_\infty$ introduces
two related difficulties.
First, the intrinsic scaling \eqref{scaling} changes
the prescribed far-field density.
Second, the evolution equation for the pressure
perturbation contains an additional term involving
the equilibrium pressure $P_\infty$, which requires
a corresponding modification of the energy estimates.

We first identify an appropriate scaling-invariant initial
quantity.
Under the intrinsic scaling
\eqref{scaling}, $D(t)$ and $P_\infty E(t)$ have
the same scaling homogeneity.
The factor $\rho_\infty^{(3-\gamma)/2}$ compensates
for this homogeneity, yielding the scaling-invariant
initial quantity $\mathcal S_0$.
For the bootstrap argument, we introduce
\begin{align}\label{xue}
\mathcal S(t)
\triangleq
\rho_\infty^{\frac{3-\gamma}{2}}
\bigl(D(t)+P_\infty E(0)\bigr),
\end{align}
which appears in
\eqref{cp1}$_2$. We note that the initial energy controls  the
pressure perturbation. More precisely, we have
\begin{align}\label{pot}
	\|p(t)\|_{L^2}^2
	\leq C_{\gamma,R}P_\infty E(0),
\end{align}
as established in \eqref{eqL}. Here,
$C_{\gamma,R}$ depends only on $\gamma$ and  $R$, rather than separately on
$\bar\rho$ and $\rho_\infty$. This estimate is a key feature
of the non-vacuum setting. Together with the basic energy
estimate \eqref{basic}, it provides the control of the
pressure perturbation needed in the subsequent higher-order
estimates.

To close the bootstrap condition \eqref{xue},
we derive an a \textit{priori} estimate for
$D(t)+P_\infty E(t)$. The key difficulty comes from the
  pressure perturbation equation:
\begin{align*}
	p_t+\operatorname{div}(pu)
	+(\gamma-1)p\operatorname{div}u
	+\gamma P_\infty\operatorname{div}u=0.
\end{align*}
The last term, which arises from the nonzero background
pressure, prevents a direct application of the arguments
developed for the far-field vacuum case in
\cite{XuZhong2026,Wen2025}. To handle this additional term,
we introduce the corrected energy functional
\begin{align}\label{hh}
	\mathscr H_K(t)
	\triangleq
	\mu\|\operatorname{curl}u(t)\|_{L^2}^2
	+\frac{1}{2\mu+\lambda}\|F(t)\|_{L^2}^2
	-\frac{2\gamma P_\infty}{2\mu+\lambda}
	\int G(\rho(t))\,dx
	+KP_\infty E(t),
\end{align}
where $K>0$ is chosen sufficiently large. Under the density
upper bound, the pressure estimate \eqref{pot}
allows us to control the negative potential-energy term in \eqref{hh}.
Consequently, for a sufficiently large choice of $K$,
\begin{align*}
\mathscr H_K(t)\sim D(t)+P_\infty E(t).
\end{align*}
Under the smallness condition \eqref{sma}, we derive
the estimate \eqref{eqh} for this corrected quantity.
By choosing $c_0$ appropriately in \eqref{c0}, we can close the bootstrap condition \eqref{xue}. The details are provided in Step 1 of Subsection \ref{sub}.

We next improve the bootstrap bound on the density, thereby closing the bootstrap argument.
Along particle trajectories, the continuity equation \eqref{1}$_{1}$
can be rewritten in the form \eqref{for}, to which
Zlotnik's inequality applies. Combining the representation of the effective viscous flux in \eqref{eqon} with the Riesz-transform estimates and the commutator estimates,  we derive the two estimates needed to apply Zlotnik's inequality; see \eqref{s1} and \eqref{s2}. Both estimates are controlled by the scaling-invariant initial quantity, see
\eqref{eqsca}. By requiring $\mathcal S_{0}$ to be
sufficiently small, Zlotnik's inequality yields the desired
uniform upper bound for the density. This establishes the desired uniform upper bound for the density and completes the bootstrap argument.

Finally, we establish the time-weighted estimates needed
to continue the local strong solution. These estimates,
together with the uniform bounds obtained above, allow us
to apply the standard continuation criterion and extend the
solution globally in time.

The remainder of the paper is organized as follows.
Section \ref{sec2} establishes the scaling-invariant
a \textit{priori} estimates, and Section \ref{sec3} derives
the time-weighted estimates.
The proof of Theorem \ref{th1} is completed in
Section \ref{sec4}.

	\section{A \textit{ priori}  estimates I: scaling invariant estimates}\label{sec2}
	For the initial data $(\rho_{0},u_{0})$ satisfies \eqref{sou},  there exists a time $T>0$ such that the problem \eqref{1}--\eqref{3} admits a unique strong solution in $\mathbb{R}^{3}\times(0,T]$, details see \cite{GLLZ20,Huang2021}. In this section, we assume that $(\rho,u)$ is a strong solution to the  problem \eqref{1}--\eqref{3} on $\mathbb R^3\times[0,T]$ for $T>0$.
 Throughout the rest of this section, $C$ denotes a generic positive constant
 which is depending only on $A, \gamma, \mu, \lambda$, and $R$.
 
 \begin{lemma}[Basic energy estimate]\label{lem1}
 	It holds that
 	\begin{align}\label{basic}
 			 \sup_{0\le t\le T}E(t)
 			+\int_0^T
 			\Big(
 			\mu\|\nabla u\|_{L^2}^{2}
 			+(\lambda+\mu)\|\operatorname{div} u\|_{L^2}^{2}
 			\Big)\, dt
 			\le CE(0).
 	\end{align}
 \end{lemma}
 
 \begin{proof}
Multiplying the momentum equation \eqref{1}$_{2}$ by $u$ and using the continuity equation \eqref{1}$_{1}$, we obtain
 	\begin{align}\label{eq1}
 		\frac{d}{d t}E(t)
 		+\mu\|\nabla u\|_{L^2}^2
 		+(\lambda+\mu)
 		\|\operatorname{div} u\|_{L^2}^2=0.
 	\end{align}
 	Integrating \eqref{eq1}	over $(0,T)$ yields \eqref{basic}.
 \end{proof}

 	\begin{lemma}\label{lem2}
 	Assume that
 	$
 	0\leq\rho\leq2\bar\rho$,
 and define
 	\begin{align*}
 		\Phi_\gamma(s)
 		\triangleq\frac{(\gamma-1)(s^\gamma-1)^2}
 		{s^\gamma-\gamma s+\gamma-1},
 		\qquad s\geq0,\quad s\neq1,
 	\end{align*}
 	with $\Phi_\gamma(1)=2\gamma$. Then $\Phi_\gamma$ is positive, continuous, and
 	nondecreasing on $[0,\infty)$, with
 	$\Phi_\gamma(0)=1$.
 	In particular,
 	\begin{align}\label{eqp}
 		P_\infty G(\rho)
 		\leq p^2
 		\leq P_\infty\Phi_\gamma(2R)G(\rho).
 	\end{align}
 	Moreover,
 	\begin{align}\label{equ}
 		\Phi_\gamma(2R)\leq C_\gamma R^\gamma,
 	\end{align}
 	where $C_\gamma>0$ depends only on $\gamma$.
 \end{lemma}
 
 \begin{proof}
 	Set
 	$
 	s\triangleq \rho/\rho_\infty$.
 	By the integral form of the mean value theorem,
 	we have
 	\begin{align}
 		p
 		=P_\infty(s^\gamma-1)
 		=P_\infty(s-1)F_\gamma(s),
 		\label{eqm}
 	\end{align}
 	where
 	\begin{align*}
 	F_\gamma(s)
 	\triangleq\gamma\int_0^1
 	\bigl(1+\theta(s-1)\bigr)^{\gamma-1}
 	\,d\theta=\frac{s^{\gamma}-1}{s-1}.
 	\end{align*}
 	Clearly,
 	$
 	F_\gamma(0)=1$, $F_\gamma(1)=\gamma$, and $F_\gamma$ is increasing on $[0,\infty)$.
 	
 	On the other hand, direct integration gives
 	\begin{align*}
 		G(\rho)
 		=\frac{A}{\gamma-1}
 		\left(
 		\rho^\gamma
 		-\gamma\rho_\infty^{\gamma-1}\rho
 		+(\gamma-1)\rho_\infty^\gamma
 		\right)=P_\infty K_\gamma(s),
 	\end{align*}
 	where
 	\begin{align*}
 	K_\gamma(s)
 	\triangleq\frac{s^\gamma-\gamma s+\gamma-1}
 	{\gamma-1}.
 	\end{align*}
 	Since
 	\begin{align*}
 	K_\gamma(1)=K_\gamma'(1)=0,
 	\qquad
 	K_\gamma''(s)=\gamma s^{\gamma-2},
 	\end{align*}
 	Taylor's formula with integral remainder
 	yields
 	\begin{align}\label{eqk}
 		K_\gamma(s)
 		=(s-1)^2 H_\gamma(s),
 	\end{align}
 	where
 	\begin{align*}
 	H_\gamma(s)
 	\triangleq\gamma\int_0^1(1-\theta)
 	\bigl(1+\theta(s-1)\bigr)^{\gamma-2}
 	\,d\theta=\frac{s^{\gamma}-\gamma s+\gamma-1}{(\gamma-1)(s-1)^{2}}.
 	\end{align*}
 	Combining \eqref{eqm} and
 	\eqref{eqk}, we obtain
 	\begin{align}\label{eqr}
 		p^2
 		=P_\infty G(\rho)
 		\frac{F_\gamma(s)^2}{H_\gamma(s)}
 		=P_\infty\Phi_\gamma(s)G(\rho).
 	\end{align}
 	Notice that $\Phi_\gamma$ extends
 	continuously to $[0,\infty)$, with
 	$
 	\Phi_\gamma(1)=2\gamma$,
 	$\Phi_\gamma(0)=1$. 
 	
 	We next prove that $\Phi_\gamma$ is
 	nondecreasing. Set
 	\begin{align}\label{d1}
 	f(s)\triangleq(s^\gamma-1)^2,
 	\qquad
 	g(s)\triangleq K_\gamma(s).
 	\end{align}
 	For $s>0$ and $s\neq1$, we have
 	\begin{align}
 		\frac{f'(s)}{g'(s)}
 		=2(\gamma-1)s^{\gamma-1}
 		\frac{s^\gamma-1}{s^{\gamma-1}-1}
 		=2\gamma s^{\gamma-1}M_\gamma(s),
 		\label{eqo}
 	\end{align}
 	where
 	\begin{align*}
 	M_{\gamma}(s)\triangleq\frac{\gamma-1}{\gamma} \frac{s^{\gamma}-1}{s^{\gamma-1}-1}.
 	\end{align*}	
 	A direct calculation gives
   $M_\gamma$ is  increasing on $(0,\infty)$.
 	Consequently, \eqref{d1} and \eqref{eqo}
 	show that
 	\begin{align*}
 	\Phi_\gamma(s)=\frac{f(s)}{g(s)}, \quad\Phi'_\gamma(s)= \frac{g^{\prime}(s)}{g(s)}\left(\frac{f^{\prime}(s)}{g^{\prime}(s)}-\frac{f(s)}{g(s)}\right)\geq 0.
 	\end{align*}
 	Since $0\leq s\leq2R$, we therefore have
 	\begin{align*}
 	1=\Phi_\gamma(0)
 	\leq\Phi_\gamma(s)
 	\leq\Phi_\gamma(2R). 
 	\end{align*}
 Combining this with \eqref{eqr}, we obtain
 \eqref{eqp}.
 	
 	Finally, we establish the bound
 	in \eqref{equ}.
 	For $s\geq2$, observe that
 	\begin{align*}
 	\frac{s^\gamma-\gamma s+\gamma-1}
 	{s^\gamma}
 	=
 	1-\gamma s^{1-\gamma}
 	+(\gamma-1)s^{-\gamma}
 	\end{align*}
 	is increasing. 	Hence
 	\begin{align*}
 	s^\gamma-\gamma s+\gamma-1
 	\geq c_\gamma s^\gamma,
 	\qquad s\geq2,
 	\end{align*}
 	where
 	\begin{align*}
 	c_\gamma
 	\triangleq1-\frac{\gamma+1}{2^\gamma}>0.
 	\end{align*}
 	It follows that
 	\begin{align*}
 	\Phi_\gamma(s)
 	\leq\frac{\gamma-1}{c_\gamma}s^\gamma,
 	\qquad s\geq2.
 	\end{align*}
 	Since $2R\geq2$, we obtain \eqref{equ}.
 \end{proof}
 
 By Lemmas \ref{lem1} and \ref{lem2}, we have the following corollary.
 \begin{corollary}
 	Under the density bound $\rho\le2\bar\rho$,
 	\begin{align}\label{eqL}
 		\|p(t)\|_{L^2}^2
 		\le\Phi_\gamma(2R)P_\infty\int G(\rho(t))\, dx
 		\le\Phi_\gamma(2R)P_\infty E(0).
 	\end{align}
 	Moreover, $P_\infty\int G\, dx$ and $\|p\|_{L^2}^2$ are equivalent,
 	with equivalence constants depending only on $R$ and $\gamma$.
 \end{corollary}

		\begin{proposition}[Bootstrap argument]
		\label{pro}
		There exist positive constants $\varepsilon_1=\varepsilon_1
		(R,A,\gamma,\mu,\lambda)$ and $c_0=c_0(R,A,\gamma,\mu,\lambda)$,
		such that the
		following statement holds. If, for all $(x,t)\in\mathbb{R}^{3}\times[0,T]$,
		\begin{align}\label{cp1}
			0\leq\rho(x,t)\leq2\bar\rho,\qquad
			\mathcal S(t)\leq2\varepsilon,
		\end{align}
		then, it holds that
		\begin{align*}
				0\leq\rho(x,t)\leq\frac32\bar\rho,\qquad
				\mathcal S(t)\leq\frac32\varepsilon,
		\end{align*}
		provided
		\begin{align}\label{cp2}
			c_0 \mathcal S_{0}\leq\varepsilon\leq\varepsilon_1.
		\end{align}
	\end{proposition}

	\begin{lemma}\label{lem3}
		Under the  assumption  \eqref{cp1}, it holds that
		\begin{align}\label{eqh}
				\sup_{0\le t\le T}\bigl(D(t)+P_\infty E(t)\bigr)
				+\int_0^T\|\sqrt\rho\dot u\|_{L^2}^2\, d t
				+P_\infty\int_0^T\|\nabla u\|_{L^2}^2\, d t
				\le C  \bigl(D(0)+P_\infty E(0)\bigr).
		\end{align}
	\end{lemma}
	\begin{proof}
		It follows from the continuity equation \eqref{1}$_{1}$
		that
		\begin{align*}
			G(\rho)_t+\operatorname{div}(G(\rho)u)
			&=G'(\rho)
			\bigl(\rho_t+u\cdot\nabla\rho\bigr)
			+G(\rho)\operatorname{div} u
			\notag\\
			&=-\bigl(\rho G'(\rho)-G(\rho)\bigr)
			\operatorname{div} u
			\notag\\
			&=-p
			\operatorname{div} u.
		\end{align*}
		Integrating the above identity
		over $\mathbb R^3$, we obtain
		\begin{align}\label{eq3}
			\frac{d}{d t}
			\int G(\rho)\, d x
			=-\int 
			p
			\operatorname{div} u \, d x.
		\end{align}

		Testing the momentum equation \eqref{1}$_{2}$ with $u_t$, then using
		$u_t=\dot u-u\cdot\nabla u$, gives
		\begin{align}\label{w0}
			&\frac12\frac d{d t}
			\left(\mu\|\nabla u\|_{L^2}^2
			+(\lambda+\mu)\|\operatorname{div} u\|_{L^2}^2\right)
			+\|\sqrt\rho\dot u\|_{L^2}^2\notag\\
			&=\int\rho\dot u\cdot(u\cdot\nabla u)\, d x
			+\frac d{d t}\int p\operatorname{div} u\, d x
			-\int p_t\operatorname{div} u\, d x.
		\end{align}
		The pressure perturbation satisfies
		\begin{align}\label{eqt}
			p_t+\operatorname{div}(pu)+(\gamma-1)p\operatorname{div} u
			+\gamma P_\infty\operatorname{div} u=0.
		\end{align}
		Multiplying \eqref{eqt} by $F$ and integrating by parts, we obtain
		\begin{align*}
		\int p_tF\, dx
		=\int p u\cdot\nabla F\, dx
		-(\gamma-1)\int p\operatorname{div} u F\, dx
		-\gamma P_\infty\int\operatorname{div} u F\, dx,
		\end{align*}
		which combine with $\operatorname{div} u=(F+p)/\nu$ implies
		\begin{align}\label{w1}
			\int p_t\operatorname{div} u\, dx
			&=\frac1{2\nu}\frac{d}{dt}\|p\|_{L^2}^2+\frac1\nu\int p_t F\, dx\notag\\
			&=\frac1{2\nu}\frac{d}{dt}\|p\|_{L^2}^2+\frac1\nu\int p\,u\cdot\nabla F\, dx
			-\frac{\gamma-1}{\nu}\int p\,\operatorname{div} u\,F\, dx
			-\frac{\gamma P_\infty}{\nu}\int \operatorname{div} u\,F\, dx.
		\end{align}		
		Denote $
			H(t)\triangleq\mu\|\operatorname{curl} u\|_{L^2}^2+\frac1\nu\|F\|_{L^2}^2$, then
		combining \eqref{w0} with \eqref{w1} yields 	
		\begin{align}
			\frac12\frac{d}{dt}H(t)+\|\sqrt\rho\dot u\|_{L^2}^2
			=&\int \rho\dot u\cdot(u\cdot\nabla u)\, dx
			-\frac1\nu\int p u\cdot\nabla F\, dx
			\notag\\
			&+\frac{\gamma-1}{\nu}\int p\operatorname{div} u\,F\, dx
			+\frac{\gamma P_\infty}{\nu}\int \operatorname{div} u\,F\, dx\notag\\
			=&\int \rho\dot u\cdot(u\cdot\nabla u)\, dx
			-\frac1\nu\int p u\cdot\nabla F\, dx
			\notag\\
			&+\frac{\gamma-1}{\nu}\int p\operatorname{div} u\,F\, dx
			+ \gamma P_\infty \|\operatorname{div}u\|_{L^2}^2
			+\frac{\gamma P_\infty}{\nu}\frac{d}{dt}\int G(\rho)\,dx,
			\label{eqH}
		\end{align}
		where we have used $\|\nabla u\|_{L^2}^2
		= \|\operatorname{curl} u\|_{L^2}^2+\|\operatorname{div} u\|_{L^2}^2$ and
		\begin{align*}
			\int \operatorname{div}u\,F\,dx
			=\nu\|\operatorname{div}u\|_{L^2}^2
			-\int p\operatorname{div}u\,dx =\nu\|\operatorname{div}u\|_{L^2}^2
			+\frac{d}{dt}\int G(\rho)\,dx.
		\end{align*}
		
 		Choose a constant $K>2\gamma/\nu$ sufficiently large, depending
 		only on the physical parameters, and recall that
 		\begin{align*}
 				\cH_K(t)&= H(t)-\frac{2\gamma P_\infty}{\nu}\int G(\rho(t))\, dx
 				+KP_\infty E(t)\\
 				&=\mu\|\operatorname{curl} u(t)\|_{L^2}^2
			+\frac1\nu\|F(t)\|_{L^2}^2
			+\frac{K P_\infty}{2}\|\sqrt\rho u(t)\|_{L^2}^2
			+\left(K-\frac{2\gamma}{\nu}\right)P_\infty \int G(\rho(t))\, dx.
		\end{align*} 
	Observe that
		\begin{align*}
		\|F\|_{L^2}^2\le C(\|\operatorname{div} u\|_{L^2}^2+\|p\|_{L^2}^2),\qquad\|p\|_{L^2}^2\le\Phi_\gamma(2R)P_\infty G(\rho),
		\qquad
		\|\operatorname{div} u\|_{L^2}^2
		\le\frac2{\nu^2}(\|F\|_{L^2}^2+\|p\|_{L^2}^2),
		\end{align*}
		there are positive constants
		$c_R$, $C_R$, depending only on $R$ and the physical parameters, such that 
		\begin{align}\label{eqeq}
				c_R\bigl[D(t)+P_\infty E(t)\bigr]
				\le\cH_K(t)
				\le C_R\bigl[D(t)+P_\infty E(t)\bigr].
		\end{align}

		Multiplying \eqref{eqH} by $2$, and adding
		$\big(KP_\infty \times \eqref{eq1}  -\frac{2\gamma}{\nu}P_\infty \times \eqref{eq3}\big)$,
		we arrive at
		\begin{align}
			&\frac{d}{dt}\cH_K(t)
			+2\|\sqrt\rho\dot u\|_{L^2}^2
			+KP_\infty\Big(
			\mu\|\nabla u\|_{L^2}^2
			+(\lambda+\mu)
			\|\operatorname{div}u\|_{L^2}^2
			\Big)\notag\\
			&=2\int 
			\rho\dot u\cdot(u\cdot\nabla u)\,dx
			-\frac{2}{\nu}\int 
			p\,u\cdot\nabla F\,dx
			+\frac{2(\gamma-1)}{\nu}
			\int 
			p\operatorname{div}u\,F\,dx
			+2\gamma P_\infty
			\|\operatorname{div}u\|_{L^2}^2.
			\label{eqe}
		\end{align}		
		We proceed to estimate the right-hand side of
		\eqref{eqe}. By the interpolation inequality and
		the assumption $\rho\leq 2\bar\rho$, we have
		\begin{align}\label{eqp6}
			\|p\|_{L^6}
			 \leq
			\|p\|_{L^\infty}^{2/3}
			\|p\|_{L^2}^{1/3}
		 \leq
			C\bar\rho^{2\gamma/3}
			\|p\|_{L^2}^{1/3}.
		\end{align}		
		Recalling the definition of the effective viscous flux,
		we take the divergence of the momentum equation
		\eqref{1}$_2$ to obtain
		\begin{align*}
			\Delta F
			=\operatorname{div}(\rho\dot u).
		\end{align*}
		 The standard elliptic and Sobolev estimates imply (details, see \cite[Lemma 2.3]{HLX2012})
		\begin{align}\label{elliptic}
		\|\nabla F\|_{L^2}\le C\|\rho\dot u\|_{L^2},\quad	\|\nabla u\|_{L^6}
			\le C\left(\bar\rho^{1/2}\|\sqrt\rho\dot u\|_{L^2}
			+\|p\|_{L^6}\right).
		\end{align}
		It follows from H\"older's, Young's inequalities, \eqref{eqp6}, and
		\eqref{elliptic}  that  
		\begin{align}
			\left|\int\rho\dot u\cdot(u\cdot\nabla u)\, dx\right|&\le \frac14\|\sqrt\rho\dot u\|_{L^2}^{2}
			+C\int\rho|u|^2|\nabla u|^2\, dx\notag\\
			&\le \frac14\|\sqrt\rho\dot u\|_{L^2}^{2}+C\bar\rho\|u\|_{L^6}^2\|\nabla u\|_{L^3}^2
			\notag\\
			&\le \frac14\|\sqrt\rho\dot u\|_{L^2}^{2}+C\bar\rho\|\nabla u\|_{L^2}^3\|\nabla u\|_{L^6}
			\notag\\
			&\le \frac14\|\sqrt\rho\dot u\|_{L^2}^{2}+C\bar\rho^{3/2}\|\nabla u\|_{L^2}^3\|\sqrt\rho\dot u\|_{L^2}
			+C\bar\rho^{1+2\gamma/3}
			\|\nabla u\|_{L^2}^3\|p\|_{L^2}^{1/3}\notag\\[1mm] 
			&\leq \frac12\|\sqrt\rho\dot u\|_{L^2}^{2}+ C \bar\rho^3D^2\|\nabla u\|_{L^2}^2
			+C\bar\rho^{1+2\gamma/3}D^{2/3}
			\|\nabla u\|_{L^2}^2.
			\label{eqp1}
		\end{align}		
		Noting that
		\begin{align*} 
			\|p\|_{L^3}^2
			\le\|p\|_{L^\infty}^{2/3}\|p\|_{L^2}^{4/3}
			\le C\bar\rho^{2\gamma/3}D^{2/3}.
		\end{align*}
		For involving $F$ terms in \eqref{eqe}, Sobolev's inequality and \eqref{elliptic}
		give
		\begin{align}\label{eqp2}
			\left|\int p\,u\cdot\nabla F\, dx\right|
			+\left|\int p\operatorname{div} u\,F\, dx\right|
			&\le C\|p\|_{L^3}\|\nabla u\|_{L^2}
			\|\nabla F\|_{L^2}\notag\\
			&\le C\bar\rho^{1/2}
			\|p\|_{L^3}\|\nabla u\|_{L^2}\|\sqrt\rho\dot u\|_{L^2}\notag\\
	     	&\le \frac12\|\sqrt\rho\dot u\|_{L^2}^{2}
			+C \bar\rho^{1+2\gamma/3}D^{2/3}
			\|\nabla u\|_{L^2}^2.	
		\end{align}
		Inserting \eqref{eqp1} and \eqref{eqp2}
		into \eqref{eqe}	yields
		\begin{align}\label{eqb1}
			&\frac{d}{dt}\cH_K(t)+\|\sqrt\rho\dot u\|_{L^2}^{2}+KP_\infty(\mu\|\nabla u\|_{L^2}^2
			+(\lambda+\mu)
			\|\operatorname{div} u\|_{L^2}^2)\notag\\
			&\le C\bar\rho^3D^2\|\nabla u\|_{L^2}^2
			+C\bar\rho^{1+2\gamma/3}D^{2/3}
			\|\nabla u\|_{L^2}^2
			+2\gamma P_\infty\|\operatorname{div} u\|_{L^2}^2.
		\end{align}
		Choose $K>1$ sufficiently large so that the last term in \eqref{eqb1}  can be absorbed, and we have
		\begin{align}\label{eqa}
		\frac{d}{d t} \mathscr{H}_{K}(t)+\|\sqrt{\rho} \dot{u}\|_{2}^{2}+  P_{\infty}\|\nabla u\|_{2}^{2} \leq C\left[\bar{\rho}^{3} D^{2}+\bar{\rho}^{1+2 \gamma / 3} D^{2 / 3}\right]\|\nabla u\|_{2}^{2}.
	    \end{align}
	    Since $\mathcal S(t)\leq 2\varepsilon$, we have
	    \begin{align*}
	    D(t)\leq
	    2\varepsilon\rho_\infty^{\frac{\gamma-3}{2}}.
	\end{align*}
	    Recalling that $\bar\rho=R\rho_\infty$ and
	    $P_\infty=A\rho_\infty^\gamma$, we obtain
	    \begin{align*}
	    	\frac{\bar\rho^3D(t)^2}{P_\infty}
	    	 \leq CR^3(2\varepsilon)^2,\quad
	    	\frac{\bar\rho^{1+\frac{2\gamma}{3}}
	    		D(t)^{\frac23}}{P_\infty}
	    	 \leq CR^{1+\frac{2\gamma}{3}}
	    	(2\varepsilon)^{\frac23}.
	    \end{align*}
One can choose  $\varepsilon_{1}$  in \eqref{cp2} sufficiently small such that
		\begin{align}\label{sma}
		C\left[R^{3}\left(2 \varepsilon_{1}\right)^{2}+R^{1+2 \gamma / 3}\left(2 \varepsilon_{1}\right)^{2 / 3}\right] \leq \frac{1}{2},
	\end{align}		
then \eqref{eqa} reduces to
		\begin{align*}
			\frac{d}{dt}\cH_K(t)+ \|\sqrt\rho\dot u\|_{L^2}^{2}
			+\frac12 P_\infty\|\nabla u(t)\|_{L^2}^2\le0.
		\end{align*}
		Integrating the above resultant over $(0,T)$ and using \eqref{eqeq} prove  \eqref{eqh}.
	\end{proof}

 \subsection{Proof of Proposition \ref{pro}}\label{sub}
 
 	\medskip
 \noindent
 \textbf{Step 1:  Closed the bootstrap estimate  for $\mathcal S(t)$.}\\
 
		By  \eqref{eqh}, we have
		\begin{align*}
			\mathcal S(t)
			&=\rho_\infty^{\frac{3-\gamma}{2}}
			\left(D(t)+P_\infty E(0)\right)\\
			&\leq
			\rho_\infty^{\frac{3-\gamma}{2}}
			\left[
			C \left(D(0)+P_\infty E(0)\right)
			+P_\infty E(0)
			\right]\\
			&\leq
			(C +1)
			\rho_\infty^{\frac{3-\gamma}{2}}
			\left(D(0)+P_\infty E(0)\right)\\[1mm]
			&=(C +1)\mathcal S_{0}.
		\end{align*}
		
		Now choose
		\begin{align}\label{c0}
			c_0\triangleq\max\left\{
			\frac12,\frac{2(C +1)}{3}
			\right\}.
		\end{align}
		Since $c_0\mathcal S_{0}\leq\varepsilon$, it follows that
		\begin{align*}
		\mathcal S(t)
		\leq(C +1)\mathcal S_{0}
		\leq\frac{C +1}{c_0}\varepsilon
		\leq\frac32\varepsilon.
		\end{align*}

	\medskip
	\noindent
	\textbf{Step 2:  Closed the bootstrap estimate  for density.}\\
	 
		For any fixed $(x,t)\in\mathbb R^3\times[0,T]$,
		let $X(s)\triangleq X(s;x,t)$ be the particle
		trajectory satisfying
		\begin{align*}
			\begin{cases}
				\displaystyle
				\frac{d}{ds}X(s)=u(X(s),s),
				&0\leq s< t,\\[1mm]
				X(t)=x.
			\end{cases}
		\end{align*}		
		For $\delta>0$, set
		\begin{align*}
				\eta_\delta(s)
				\triangleq\delta\exp\left(
				-\int_0^s\operatorname{div}u(X(\tau),\tau)
				\,d\tau\right),\quad
				\rho_\delta(s)
				 \triangleq\rho(X(s),s)+\eta_\delta(s)>0.
		\end{align*}
		By the continuity equation \eqref{1}$_{1}$,
		\begin{align*}
		\frac{d}{ds}\rho(X(s),s)
		=-\rho(X(s),s)\operatorname{div}u(X(s),s),
	\end{align*}
		while
		$
		\eta_\delta'(s)
		=-\eta_\delta(s)\operatorname{div}u(X(s),s)$.
		Therefore,
		\begin{align*}
			\frac{d}{ds}\rho_\delta(s)
			+\rho_\delta(s)\operatorname{div}u(X(s),s)=0.
		\end{align*}
		This gives
		\begin{align}\label{for}
		Y_\delta'(s)=g_\delta(s)+b^{\prime}(s),
	\end{align}
		where
	 \begin{align*}
			 Y_\delta(s)\triangleq\log\rho_\delta(s),\quad b(s)\triangleq-\frac1\nu\int_0^s
			F(X(\tau),\tau)\,d\tau,\quad
			g_\delta(s)\triangleq-\frac{A\rho^{\gamma}(X(s),s) -P_{\infty}}{\nu}.
		\end{align*}

		Define
		\begin{align*}
	      \phi\triangleq(-\Delta)^{-1}\operatorname{div}(\rho u).
		\end{align*}	 
		Since
		$
		\rho\dot u
		=(\rho u)_t+\operatorname{div}(\rho u\otimes u)$
		and
		$
		\Delta F=\operatorname{div}(\rho\dot u)$,
		we obtain
		\begin{align*}
			F
		=-\phi_t
			-\mathcal R_{ij}(\rho u^iu^j)=-\bigl(\phi_t+u\cdot\nabla\phi\bigr)
			+[u^i,\mathcal R_{ij}](\rho u^j),
		\end{align*}
		where
		\begin{align*}
		\mathcal R_{ij}
		\triangleq\partial_i(-\Delta)^{-1}\partial_j,\quad
		[u^i,\mathcal R_{ij}](\rho u^j)
		\triangleq u^i\mathcal R_{ij}(\rho u^j)
		-\mathcal R_{ij}(\rho u^iu^j).
	\end{align*}
Along the trajectory $X(s)$,
		\begin{align}\label{eqon}
			F(X(s),s)
			=-\frac{d}{ds}\phi(X(s),s)
			+[u^i,\mathcal R_{ij}](\rho u^j)(X(s),s).
		\end{align}
		It follows that, for
		$0\leq t_1<t_2\leq t$,
		\begin{align}
			b(t_2)-b(t_1)
			 \leq\frac2\nu
			\sup_{0\leq s\leq T}
			\|\phi(s)\|_{L^\infty}
		 +\frac1\nu\int_{t_1}^{t_2}
			\|[u^i,\mathcal R_{ij}](\rho u^j)\|_{L^\infty}
			\,ds
			 \triangleq J_1+J_2.
			\label{eqb}
		\end{align}
		
		By \eqref{basic}, the Gagliardo--Nirenberg inequality,
		and the boundedness of the Riesz transforms on
		$L^q(\mathbb R^3)$ for $1<q<\infty$, we obtain
		\begin{align*}
			\|\phi\|_{L^\infty}
		 \leq
			C\|\nabla\phi\|_{L^2}^{1/2}
			\|\nabla\phi\|_{L^6}^{1/2} \leq
			C\|\rho u\|_{L^2}^{1/2}
			\|\rho u\|_{L^6}^{1/2} \leq
			C\bar\rho^{3/4}E(0)^{1/4}
			\|\nabla u\|_{L^2}^{1/2}.
		\end{align*}
		Consequently, it follows from \eqref{eqh} that
		\begin{align}\label{eqJ1}
			J_1
			\leq
			C\left[
			\bar\rho^3 E(0)
			\bigl(D(0)+P_\infty E(0)\bigr)
			\right]^{1/4}.
		\end{align}
		Applying the commutator estimates (details, see Proposition 2.6 in \cite{Li2020}), Cauchy--Schwarz's inequality, \eqref{eqh}, and	\eqref{eqp6}, we have
		\begin{align}
			J_2
			&\leq C\bar\rho\int_{t_1}^{t_2}
			\|\nabla u\|_{L^2}
			\|\nabla u\|_{L^6}\,ds
			\notag\\
			&\leq C\bar\rho^{3/2}
			\int_{t_1}^{t_2}
			\|\nabla u\|_{L^2}
			\|\sqrt\rho\dot u\|_{L^2}\,ds
			+C\bar\rho
			\int_{t_1}^{t_2}
			\|\nabla u\|_{L^2}\|p\|_{L^6}\,ds
			\notag\\
			 &\leq C\bar\rho^{3/2}
			 \left(\int_0^T\|\nabla u\|_{L^2}^2\,ds
			 \right)^{1/2}
			 \left(\int_0^T
			 \|\sqrt\rho\dot u\|_{L^2}^2\,ds
			 \right)^{1/2}+C\bar\rho^{1+ 2\gamma/3}
			 \int_{t_1}^{t_2}
			 \|\nabla u\|_{L^2}
			 \|p\|_{L^2}^{1/3}\,ds\notag\\[2mm] 
			 &\leq C
			 \bigl[\bar\rho^3E(0)\bigl(D(0)+P_\infty E(0)\bigr)\bigr]^{1/2}+C\bar\rho^{1+ 2\gamma/3}
			 \int_{t_1}^{t_2}D(s)^{2/3}\,ds\notag\\[2mm] 
			 &\leq C
			 \bigl[\bar\rho^3E(0)\bigl(D(0)+P_\infty E(0)\bigr)\bigr]^{1/2}+C\bar\rho^{1+ 2\gamma/3}
			 \bigl(D(0)+P_\infty E(0)\bigr)^{2/3}(t_2-t_1).\label{eqJ2}
		\end{align}

		It follows from \eqref{eqb}--\eqref{eqJ2} that
		\begin{align*}
			b(t_2)-b(t_1)
			\leq N_0+N_1(t_2-t_1),
		\end{align*}
		where
		\begin{gather} 
			N_1
			\triangleq C\bar\rho^{1+ 2\gamma/3}
			\bigl(D(0)+P_\infty E(0)\bigr)^{2/3},\label{s1}\\[2mm] 
			N_0
			\triangleq C \left[
			\bigl[\bar\rho^3E(0)\bigl(D(0)+P_\infty E(0)\bigr)\bigr]^{1/4}
			+\bigl[\bar\rho^3E(0)\bigl(D(0)+P_\infty E(0)\bigr)\bigr]^{1/2}
			\right].\label{s2}
		\end{gather}
		
	The definition of $\mathcal S_{0}$ implies
		\begin{align*}
		\rho_\infty^{\frac{\gamma+3}{2}}E(0)
		\leq A^{-1}\mathcal S_{0}.
		\end{align*}
		Then,
		\begin{align*}
			\bar\rho^3E(0)\bigl(D(0)+P_\infty E(0)\bigr)
			 =R^3
			\left(\rho_\infty^{\frac{\gamma+3}{2}}E(0)\right)
			\left[\rho_\infty^{\frac{3-\gamma}{2}}\bigl(D(0)+P_\infty E(0)\bigr)\right] 
			 \leq  C R^3\mathcal S_{0}^2,
		\end{align*}
		and
		\begin{align*}
			\bar\rho^{1+ 2\gamma/3}\bigl(D(0)+P_\infty E(0)\bigr)^{2/3}
			 =
			R^{1+ 2\gamma/3}
			\rho_\infty^\gamma \mathcal S_{0}^{2/3} =
			A^{-1}P_\infty
			R^{1+ 2\gamma/3}\mathcal S_{0}^{2/3}.
		\end{align*}
		Thus,  
		\begin{align}\label{eqsca}
			N_0
			 \leq C \left(
			 \mathcal S_{0}^{1/2}
			+ \mathcal S_{0}
			\right),
			\quad
			N_1
			&\leq C P_\infty \mathcal S_{0}^{2/3}.
		\end{align}
		
		We next apply Zlotnik's inequality (see Lemma \ref{lem:zlotnik}) to obtain
		the desired uniform density bound.
		Recalling the definition of $g_\delta(s)$ that
		\begin{align*}
			g_\delta(Y_\delta)
			=-\frac1\nu
			\left[
			A\bigl(e^{Y_\delta}-\eta_\delta(s)\bigr)^\gamma
			-P_\infty
			\right].
		\end{align*}
		Moreover,
		$
		\lim_{y\to+\infty}g_\delta(y)=-\infty$. Choose
		\begin{align*}
			\zeta_\delta
			\triangleq\log\left(
			\frac54\bar\rho+\delta \exp\Big(
			\int_0^T
			\|\operatorname{div}u(X(\tau),\tau)\|_{L^{\infty}}\,d\tau
			\Big)
			\right).
		\end{align*}
		For any $Y_\delta\geq\zeta_\delta$, using
		$\rho_\infty\leq\bar\rho$, we have
		\begin{align*}
			g_\delta(Y_\delta)
			\leq
			-\frac{A}{\nu}
			\left[
			\left(\frac54\bar\rho\right)^\gamma
			-\rho_\infty^\gamma
			\right]
		\leq
			-\frac{A}{\nu}
			\left[
			\left(\frac54\right)^\gamma-1
			\right]\bar\rho^\gamma.
		\end{align*}
		By \eqref{eqsca}, we  choose
		$\varepsilon_{1}$ in \eqref{cp2}  sufficiently small such that
		\begin{align}\label{eqZ}
			N_0\leq\log\frac65,
			\qquad
			N_1\leq
			\frac{A}{\nu}
			\left[
			\left(\frac54\right)^\gamma-1
			\right]\bar\rho^\gamma.
		\end{align}
		In particular,
		$
		g_\delta(Y_\delta)\leq-N_1$,
		 for all $Y_\delta\geq\zeta_\delta$.\\[-0.5mm]

		Combining \eqref{for} and
		\eqref{eqZ}, we apply
		the Zlotnik's inequality to obtain
		\begin{align*} 
			Y_\delta(s)
			&\leq
			\max\{Y_\delta(0),\zeta_\delta\}+N_0\\[2mm] 
			&\leq
			\max\left\{
			\log(\bar\rho+\delta),
			\zeta_\delta
			\right\}
			+\log\frac65.
		\end{align*}
		Letting $\delta\to0^+$, we obtain
		\begin{align*}
		\rho(x,t)\leq\frac32\bar\rho.
	    \end{align*}
		This completes the proof. \hfill$\Large\Box$
	
\section{A \textit{ priori}  estimates II: time weighted estimates}\label{sec3}

Throughout this subsection, we assume that $(\rho,u)$ is a strong solution to the problem \eqref{1}--\eqref{3} in $\mathbb{R}^{3}\times(0,T]$ for some $T>0$, and the bootstrap assumptions
\eqref{cp1} hold.  As a direct consequence of the estimates obtained in Section \ref{sec2}, we have
\begin{align}\label{you}
	0\le\rho(x,t)\le 2\bar\rho,
	\qquad
	\sup_{0\le t\le T}\bigl(E(t)+D(t)\bigr)
	+\int_0^T\bigl(\|\nabla u\|_2^2+\|\sqrt\rho\dot u\|_2^2\bigr)\, dt
	\le C,
\end{align}
where $C$ denotes a generic positive constant depending on $\bar\rho, \rho_\infty, A, \gamma, \mu, \lambda, q$, and the initial data. We will derive higher-order estimates which are needed to establish the time continuity of the solution.
The letter $C(T)$ may depend on a fixed finite time $T$.

\begin{proposition}\label{prop}
For  finite $T<\infty$, it holds that
	\begin{align*}
		\sup_{0<t\le T}\left( t^i\|\sqrt\rho\dot u(t)\|_{L^2}^2\right) 
		+\int_0^T t^i\|\nabla\dot u\|_{L^2}^2\, dt\le C(T),
		\qquad i=1,2.
	\end{align*}
	Moreover,
	\begin{align*}
		&\sup_{0\le t\le T}\left(\|\rho-\rho_\infty\|_{H^1\cap W^{1,q}}+
		\|\rho_t\|_{L^2}\right)
		+\int_0^T\left(\|\nabla u\|_{L^{\infty}}+
		\|\nabla^2u\|_{L^2}+\|\nabla^2u\|_{L^q}+\|\rho\dot u\|_{L^q}\right)\, dt\le C(T),\\
		&\sup_{0<t\le T}\left( t\|\nabla^2u(t)\|_{L^2}^2\right) 
		+\int_0^Tt\left(\|\nabla u_t\|_{L^2}^2+
		\|\nabla^2u\|_{L^q}^2+\|\rho\dot u\|_{L^q}^2\right)\, dt\le C(T).
	\end{align*}
\end{proposition}
\begin{proof}
The proof is standard and follows from the time-weighted
energy estimates, the Lam\'e elliptic estimates, and
the Beale-Kato-Majda-type inequality.
For brevity, we omit the details and refer the reader
to \cite{HLX2012,XuZhong2026} for similar arguments.
\end{proof}

\section{Proof of Theorem \ref{th1}}\label{sec4}
We are now in a position to prove Theorem \ref{th1}. By the local
well-posedness result have been proved in \cite{GLLZ20,Huang2021}, there exists a time
$T_{*}>0$ such that the problem
\eqref{1}--\eqref{3} admits a unique strong solution $(\rho,u)$ on
$\mathbb{R}^{3}\times(0,T_{*}]$. We first verify that the bootstrap assumptions are valid for a short time.
Indeed, by \eqref{sou}, \eqref{sca}, and the choice of $c_0$, the initial data satisfy \eqref{cp1} at $t=0$.
Therefore, by continuity, there exists $T_1\in(0,T_*]$ such that
\eqref{cp1} holds on $[0,T_1]$.

Define
\begin{align*}
	T^{*}\triangleq \sup \left\{ T>0 \,\middle|\,
	(\rho, u) \text{ is a strong solution of } \eqref{1}-\eqref{3}
	\text{ on } \mathbb{R}^{3}\times[0,T] \text{ satisfying } \eqref{cp1}
	\right\}.
\end{align*}
Clearly, $T^{*}\geq T_1>0$.   For any $0<\tau<T\leq T^{*}$ with $T<\infty$, Proposition \ref{prop} and \eqref{you} give
\begin{align*}
	\rho u\in C\bigl([0,T];L^2\bigr),
	\qquad
	\rho-\rho_{\infty}\in C\bigl([0,T];H^1\cap W^{1,q}\bigr).
\end{align*}
If $T^{*}<\infty$, the limit
\begin{align*}
	(\rho, \rho u)(T^{*})
	=
	\lim_{t\to T^{*}}(\rho, \rho u)(t)
\end{align*}
holds. The limiting data at time $T^{*}$ also satisfy the  condition
required in \eqref{sou}.
Taking
$(\rho, \rho u)(T^{*})$ as a new initial data, then the local
well-posedness result yields
a strong solution on $[T^{*},T^{*}+\delta]$ for some $\delta>0$.  This contradicts
the definition of $T^{*}$. Therefore, $T^{*}=\infty$, and the strong solution exists globally in
time. The uniqueness follows from the local uniqueness and the continuation
argument above. This completes the proof of Theorem \ref{th1}. \hfill$\square$

\section{Appendix}	
 
\setcounter{lemma}{0}
\renewcommand{\thelemma}{A.\arabic{lemma}}	 
\begin{lemma}\label{check}
The smallness condition
\begin{align}\label{ch1}
	\mathcal S_0
	\leq
	\varepsilon_0\left(
	\frac{\bar\rho}{\rho_\infty}
	\right)
\end{align}
is invariant under the scaling transformation
\eqref{scaling}.
\end{lemma}

	\begin{proof}
	
		The far-field density and pressure satisfy
		\begin{align}\label{eqr1}
			\rho_\infty^\ell
			=
			\ell^{\frac{2}{\gamma+1}}\rho_\infty, \quad
				P_\infty^\ell
				 =
				A(\rho_\infty^\ell)^\gamma =
				A\left(
				\ell^{\frac{2}{\gamma+1}}\rho_\infty
				\right)^\gamma =
				\ell^{\frac{2\gamma}{\gamma+1}}
				P_\infty.
		\end{align}	
		For convenience, we denote
		\begin{align*}
		t_\ell
		\triangleq
		\ell^{\frac{2\gamma}{\gamma+1}}t, \quad y\triangleq\ell x.
		\end{align*}
		By \eqref{scaling}, we have
		\begin{align}\label{gen}
			\frac12
			\int 
			\rho^\ell(x,t)|u^\ell(x,t)|^2\,dx
			=
			\frac12
			\ell^{\frac{2\gamma}{\gamma+1}}
			\int 
			\rho(\ell x,t_\ell)
			|u(\ell x,t_\ell)|^2\,dx=
			\frac12
				\ell^{-\frac{\gamma+3}{\gamma+1}}
			\int 
			\rho(y,t_\ell)|u(y,t_\ell)|^2\,dy.
		\end{align}
	The potential corresponding to the scaled density is defined by
		\begin{align*}
		G^\ell(\rho^\ell)
		=
		\rho^\ell
		\int_{\rho_\infty^\ell}^{\rho^\ell}
		\frac{P(s)-P(\rho_\infty^\ell)}{s^2}\,ds.
		\end{align*}
		Set
		\begin{align*}
		s\triangleq\ell^{\frac{2}{\gamma+1}}r.
		\end{align*}
		Then
		\begin{align*}
			P(s)-P(\rho_\infty^\ell)
			=
			A\left[
			\ell^{\frac{2\gamma}{\gamma+1}}r^\gamma
			-
			\ell^{\frac{2\gamma}{\gamma+1}}
			\rho_\infty^\gamma
			\right]=
			\ell^{\frac{2\gamma}{\gamma+1}}
			\left[P(r)-P(\rho_\infty)\right].
		\end{align*}
		Therefore,
		\begin{align*}
			G^\ell(\rho^\ell)
			=
			\ell^{\frac{2}{\gamma+1}}\rho
			\int_{\rho_\infty}^{\rho}
			\frac{
				\ell^{\frac{2\gamma}{\gamma+1}}
				[P(r)-P(\rho_\infty)]
			}{
				\ell^{\frac{4}{\gamma+1}}r^2
			}
			\ell^{\frac{2}{\gamma+1}}\,dr=
			\ell^{\frac{2\gamma}{\gamma+1}}
			\rho
			\int_{\rho_\infty}^{\rho}
			\frac{P(r)-P(\rho_\infty)}{r^2}\,dr=
			\ell^{\frac{2\gamma}{\gamma+1}}
			G(\rho).
		\end{align*}
		Consequently,
		\begin{align}\label{gen1}
			\int 
			G^\ell(\rho^\ell(x,t))\,dx
		=
			\ell^{\frac{2\gamma}{\gamma+1}}
			\int 
			G(\rho(\ell x,t_\ell))\,dx=
			\ell^{\frac{2\gamma}{\gamma+1}-3}
			\int 
			G(\rho(y,t_\ell))\,dy=
					\ell^{-\frac{\gamma+3}{\gamma+1}}
					\int 
					G(\rho(y,t_\ell))\,dy.
				\end{align}
		Combining \eqref{gen} and
		\eqref{gen1}, we deduce that
		\begin{align}\label{eqE}
				E^\ell(0)
				=
				\ell^{-\frac{\gamma+3}{\gamma+1}}E(0).
		\end{align}

		Differentiating the scaled velocity with respect
		to the spatial variable gives
		\begin{align}\label{eqg}
			\|\nabla u^\ell(t)\|_{L^2}^2
			=
			\ell^{\frac{4\gamma}{\gamma+1}}
			\int 
			|\nabla u(\ell x,t_\ell)|^2\,dx=
			\ell^{\frac{\gamma-3}{\gamma+1}}
			\int 
			|\nabla u(y,t_\ell)|^2\,dy
				=
				\ell^{\frac{\gamma-3}{\gamma+1}}
				\|\nabla u(t_\ell)\|_{L^2}^2.
		\end{align}		
		Recall that
		$
		p=P(\rho)-P_\infty$.
		For the scaled solution, we define
		$
		p^\ell
	\triangleq
		P(\rho^\ell)-P_\infty^\ell$.
		By \eqref{scaling}, we have
		\begin{align*}
			p^\ell(x,t)
			=
			A(\rho^\ell(x,t))^\gamma
			-
			P_\infty^\ell=
			A\ell^{\frac{2\gamma}{\gamma+1}}
			\rho^\gamma(\ell x,t_\ell)
			-
			\ell^{\frac{2\gamma}{\gamma+1}}P_\infty=
			\ell^{\frac{2\gamma}{\gamma+1}}
		\left( 
			P(\rho(\ell x,t_\ell))-P_\infty
			\right) =
			\ell^{\frac{2\gamma}{\gamma+1}}
			p(\ell x,t_\ell).
		\end{align*}		
		It follows that
		\begin{align}\label{eqpp}
			\|p^\ell(t)\|_{L^2}^2
			=
			\ell^{\frac{4\gamma}{\gamma+1}}
			\int 
			|p(\ell x,t_\ell)|^2\,dx=
			\ell^{\frac{4\gamma}{\gamma+1}-3}
			\int 
			|p(y,t_\ell)|^2\,dy
				=
				\ell^{\frac{\gamma-3}{\gamma+1}}
				\|p(t_\ell)\|_{L^2}^2.
		\end{align}	
		Combining \eqref{eqg} and
		\eqref{eqpp}, we obtain
		\begin{align}\label{DD}
			D^\ell(t)
		=
			\|\nabla u^\ell(t)\|_{L^2}^2
			+\|p^\ell(t)\|_{L^2}^2=
			\ell^{\frac{\gamma-3}{\gamma+1}}
			\left(
			\|\nabla u(t_\ell)\|_{L^2}^2
			+\|p(t_\ell)\|_{L^2}^2
			\right)=
			\ell^{\frac{\gamma-3}{\gamma+1}}
			D(t_\ell).
		\end{align}

		By \eqref{eqr1} and
		\eqref{eqE}, we have
			\begin{align}\label{eqP}
			P_\infty^\ell E^\ell(0)
			=
			\ell^{\frac{2\gamma}{\gamma+1}}
			P_\infty
			\ell^{-\frac{\gamma+3}{\gamma+1}}
			E(0)=
			\ell^{
				\frac{2\gamma-(\gamma+3)}{\gamma+1}
			}
			P_\infty E(0)=
			\ell^{\frac{\gamma-3}{\gamma+1}}
			P_\infty E(0).
		\end{align}		
		It follow from \eqref{eqE}, \eqref{DD}, \eqref{eqP} that
		\begin{align}\label{ss2}
				\mathcal S_{0}^\ell =\mathcal S_{0},
				\qquad\text{for every }\ell>0.
		\end{align}	
		Hence, $\mathcal S_0$ is invariant under the scaling
		transformation \eqref{scaling}.	It follows from \eqref{ss} and \eqref{ss2} that lemma holds.
		\end{proof}

		\setcounter{lemma}{1}
		\renewcommand{\thelemma}{A.\arabic{lemma}}	 
		\begin{lemma}\label{lem:zlotnik}\cite[Lemma 1.3]{Z2000}
			Let the function  $y$  satisfy
			\begin{align*}
				y^{\prime}(t)=g(y)+b^{\prime}(t) \text { on }[0, T],\quad  y(0)=y^{0},
			\end{align*}
			with  $g \in C(R)$  and  $y$, $b \in W^{1,1}(0, T)$. If  $g(\infty)=-\infty$  and
			\begin{align*}
				b\left(t_{2}\right)-b\left(t_{1}\right) \leq N_{0}+N_{1}\left(t_{2}-t_{1}\right),
			\end{align*}
			for all  $0 \leq t_{1}<t_{2} \leq T$  with some  $N_{0} \geq 0$  and  $N_{1} \geq 0$, then
			\begin{align*}
				y(t) \leq \max \left\{y^{0}, \bar{\zeta}\right\}+N_{0}<\infty\quad \text { on }\quad [0, T],
			\end{align*}
			where  $\bar{\zeta}$  is a constant such that
			\begin{align*}
				g(\zeta) \leq-N_{1} \quad\text { for }\quad \zeta \geq \bar{\zeta}.
			\end{align*}
		\end{lemma}

	\section*{Conflict of interest}
	The authors declare that they have no conflict of interest.
	
	\section*{Data availability}
	No data was used for the research described in the article.

	\section*{AI declaration}
	The authors declare that no AI tools have been used in the production of the manuscript. All contents, analysis,  writing, and reasoning are our original independent work.
	
	


\begin{thebibliography}{99}
		
		\bibitem{CD2010}
		F. Charve  and R. Danchin, A global existence result for the compressible Navier--Stokes equations in the critical $L^{p}$ framework, {\it Arch. Ration. Mech. Anal.}, {\bf 198} (2010), 233--271.
		
		\bibitem{CMZ2010}
		Q. Chen,  C. Miao,  and Z. Zhang, Global well-posedness for compressible Navier--Stokes equations with highly oscillating initial velocity, {\it Comm. Pure Appl. Math.}, {\bf 63} (2010), 1173--1224.
		
		\bibitem{Danchin2000}
		R. Danchin, Global existence in critical spaces for compressible Navier--Stokes equations, {\it Invent. Math.}, {\bf 141} (2000), 579--614.
		
		\bibitem{FZZ2018}
		D. Fang, T. Zhang, and R. Zi, Global solutions to the isentropic compressible Navier--Stokes equations with a class of large initial data, {\it SIAM J. Math. Anal.}, {\bf 50} (2018), 4983--5026.
		
		\bibitem{GLLZ20}
		H. Gong,  J. Li, X. Liu, and X. Zhang, Local well-posedness of isentropic compressible Navier--Stokes equations with vacuum, {\it Commun. Math. Sci.},  {\bf 18} (2020), 1891--1909.
		
	\bibitem{H2011a}
	B. Haspot, Well-posedness in critical spaces for the system of compressible Navier--Stokes in larger spaces, {\it J. Differential Equations}, {\bf 251} (2011), 2262--2295.
	
	\bibitem{H2011b}
	B. Haspot, Existence of global strong solutions in critical spaces for barotropic viscous fluids, {\it Arch. Ration. Mech. Anal.}, {\bf 202} (2011), 427--460.
	
	
	\bibitem{HHW2019}
L. He, J. Huang, and C. Wang, Global stability of large solutions to the 3D compressible Navier--Stokes equations,  {\it Arch. Ration. Mech. Anal.}, {\bf 234} (2019), 1167--1222.
	
	 
	
	\bibitem{HHPZ2024}
	G. Hong, X. Hou, H. Peng, and C. Zhu, Global existence for a class of large solution to compressible Navier--Stokes equations with vacuum,  {\it Math. Ann.}, {\bf 388}   (2024), 2163--2194.
		
	\bibitem{Huang2021}
	X. Huang,  On local strong and classical solutions to the three-dimensional barotropic compressible Navier--Stokes equations with vacuum, {\it Sci. China Math.}, {\bf 64} (2021), 1771--1788.
	
	\bibitem{HL2022}
	X. Huang  and J. Li, Global well-posedness of classical solutions to the Cauchy problem of two-dimensional barotropic compressible Navier--Stokes system with vacuum and large initial data, {\it SIAM J. Math. Anal.}, {\bf 54} (2022), 3192--3214.

	

	
	\bibitem{HLX2012}
	X. Huang, J. Li, and Z. Xin, Global well-posedness of classical solutions with large oscillations and vacuum to the three-dimensional isentropic compressible Navier--Stokes equations, {\it Comm. Pure Appl. Math.}, {\bf 65} (2012), 549--585.
	
	\bibitem{LX2019}
	Z. Lei and Z. Xin, On scaling invariance and type-I singularities for the compressible Navier--Stokes equations, {\it Sci. China  Math.}, {\bf 62} (2019), 2271--2286.
	
	\bibitem{Li2020}
	J. Li, Global small solutions of heat conductive compressible Navier--Stokes equations with vacuum: smallness on scaling invariant quantity, {\it Arch. Ration. Mech. Anal.}, {\bf 237} (2020), 899--919.
	
	\bibitem{LiXin2019}
	J. Li and Z. Xin, Global well-posedness and large time asymptotic behavior of classical solutions to the compressible Navier--Stokes equations with vacuum, {\it Ann. PDE}, {\bf 5} (2019), Paper No. 7.
	
	
	\bibitem{MN1980}
	A. Matsumura and T. Nishida, The initial value problem for the equations of motion of viscous and heat-conductive gases, {\it J. Math. Kyoto Univ.}, {\bf 20} (1980), 67--104.
	

	\bibitem{WWZ2014}
	C. Wang, W. Wang, and Z. Zhang, Global well-posedness of compressible Navier--Stokes equations for some classes of large initial data, {\it Arch. Ration. Mech. Anal.}, {\bf 213} (2014), 171--214.
	
	\bibitem{Wen2025}
	H. Wen, Global wellposedness of compressible Navier--Stokes equations with vacuum and smallness on scaling invariant quantity in $\mathbb R^3$, {\it Adv. Math.},  {\bf 482} (2025), Paper No. 110628.
	
	\bibitem{XXZ2026}
	M. Xie, S. Xu, and Y. Zhang, Global strong solutions of compressible Navier--Stokes equations with small scaling invariant, {\it J. Math. Phys.}, {\bf 67} (2026), Paper No. 081504.
	
	\bibitem{XZ20261}
	H. Xu and X. Zhong, Global well-posedness for the compressible Navier--Stokes equations with vacuum and smallness on coefficient-coupled scaling invariant quantity, {\it J. Dynam. Differential Equations}, doi: 10.1007/s10884-026-10548-5, 2026.
	
	\bibitem{XuZhong2026}
	L. Xu and X. Zhong, Pressure effects on critical scaling and global low-regularity solutions for compressible Navier--Stokes system, https://arxiv.org/abs/2606.11575.
	 
	\bibitem{XZ202611}
	L. Xu and X. Zhong, Global well-posedness of isentropic compressible Navier--Stokes equations with smallness on scaling-invariant quantity in a half-space, https://arxiv.org/abs/2607.21926.
	

     
     \bibitem{ZLZ2020}
     X. Zhai, Y. Li, and F. Zhou, Global large solutions to the three dimensional compressible Navier--Stokes equations, {\it SIAM J. Math. Anal.}, {\bf 52} (2020), 1806--1843.
     
     \bibitem{Z2000}
     A. A. Zlotnik, Uniform estimates and the stabilization of symmetric solutions of a system of quasilinear equations,  {\it Differ. Equ.}, {\bf 36} (2000),  701--716.
	
	
\end{thebibliography}
\end{document}